\documentclass[a4paper,twoside,11pt]{article}
\usepackage{bbm}
\usepackage{amsfonts}
\usepackage{amssymb}\usepackage{amsmath}
 \numberwithin{equation}{section}
\begin{document}

\title{Global weak solutions to a weakly dissipative $\mu$HS equation}

\author{Jingjing Liu\footnote{e-mail:
jingjing830306@163.com}\ \ \ \ \ \ \ \ Zhaoyang Yin\footnote{e-mail: mcsyzy@mail.sysu.edu.cn}
\\ Department of Mathematics, Sun Yat-sen University,\\ 510275 Guangzhou, China}
\date{}
\maketitle

\begin{abstract}
This paper is concerned with global existence of weak solutions
 for a weakly dissipative $\mu$HS equation by using smooth
approximate to initial data and Helly$^{,}$s theorem. \\
\end{abstract}


\section{Introduction}
\newtheorem {definition1}{Definition}[section]
\newtheorem{theorem1}{Theorem}[section]
\newtheorem{remark1}{Remark}[section]
\par
Recently, the $\mu$-Hunter-Saxton (also called $\mu$-Camassa-Holm) equation
$$\mu(u)_{t}-u_{txx}=-2\mu(u)u_{x}+2u_{x}u_{xx}+uu_{xxx},$$
which is originally derived and studied in \cite{k-l-m} attracts a lot of attention. Here
$u(t,x)$ is a time-dependent function on the unit circle $\mathbb{S}=\mathbb{R}/\mathbb{Z}$ and
$\mu(u)=\int_{\mathbb{S}}udx$ denotes its mean. The
$\mu$-Hunter-Saxton equation lies mid-way between the periodic
Hunter-Saxton and Camassa-Holm equations. In \cite{k-l-m}, the authors show that if interactions
of rotators and an external magnetic field is allowed, then the $\mu$HS equation can
be viewed as a natural generalization of the rotator equation. Moreover, the $\mu$HS equation describes the geodesic flow
on $\mathcal{D}^{s}(\mathbb{S})$ with the right-invariant metric given at the
identity by the inner product \cite{k-l-m}
$$(u,v)=\mu(u)\mu(v)+\int_{\mathbb{S}}u_{x}v_{x}dx.$$ In \cite{k-l-m, l-m-t}, the authors showed that the $\mu$HS equation admits both periodic one-peakon solution and the multi-peakons. Moreover, in \cite{fu, gui}, the authors also discussed the $\mu$HS equation.

 In general, it is difficult to avoid energy dissipation mechanisms in a real world. So, it is reasonable to study the model with energy dissipation. In \cite{g} and \cite{os}, the authors discussed the energy dissipative KdV equation from different aspects. Weakly dissipative CH equation and weakly dissipative DP equation have been studied in \cite{wu, wy1} and \cite{ewy, g-l-w, wy, y1} respectively. Recently, some results for a weakly dissipative $\mu$DP equation are proved in \cite{k2}. It is worthy to note that the local well-posedness result in \cite{k2} is obtained by using a method based on a geometric argument.

 In this paper, we will discuss global existence of weak solutions of the following weakly dissipative $\mu$HS equation:
\begin{equation}
\left\{\begin{array}{ll}
y_{t}+uy_{x}+2u_{x}y+\lambda y=0,&t > 0,\,x\in \mathbb{R},\\
y=\mu(u)-u_{xx},&t > 0,\,x\in \mathbb{R},\\
u(0,x) = u_{0}(x),& x\in \mathbb{R}, \\
u(t,x+1)=u(t,x), & t \geq 0, x\in \mathbb{R},\\\end{array}\right. \\
\end{equation}
or in the equivalent form:
\begin{equation}
\left\{\begin{array}{ll}
\mu(u)_{t}-u_{txx}+2\mu(u)u_{x}-2u_{x}u_{xx}-uu_{xxx}+\lambda(\mu(u)-u_{xx})=0,&t > 0,\,x\in \mathbb{R},\\
u(0,x) = u_{0}(x),& x\in \mathbb{R}, \\
u(t,x+1)=u(t,x), & t \geq 0, x\in \mathbb{R}.\\\end{array}\right. \\
\end{equation}
Here the constant $\lambda$ is assumed to be
positive and the term $\lambda y=\lambda(\mu(u)-u_{xx})$ models energy dissipation. The Cauchy problem of (1.1) has been discussed in
\cite{l-y} recently. The authors established the local
well-posedness, derived the precise blow-up scenario for Eq.(1.1) and proved that Eq.(1.1) has global strong solutions
and also finite time blow-up solutions. However, the existence of
global weak solutions to Eq.(1.1) has not been studied yet.
The aim of this paper is to present a  global existence result of
weak solutions to Eq.(1.1).

Throughout the paper, we denote by $*$ the convolution. Let $\|\cdot\|_{Z}$ denote the norm of Banach space $Z$ and let $\langle \cdot, \cdot\rangle$
denote the $H^{1}(\mathbb{S})$, $H^{-1}(\mathbb{S})$ duality bracket. Let $M(\mathbb{S})$ be the space of Radon
measures on $\mathbb{S}$ with bounded total variation and $M^{+}(\mathbb{S})$ ($M^{-}(\mathbb{S})$) be the subset of $M(\mathbb{S})$ with positive (negative) measures. Finally, we write $BV(\mathbb{S})$ for the space of functions with bounded
variation, $\mathbb{V}(f)$ being the total variation of $f\in BV(\mathbb{S}).$

Before giving the precise statement of the main result, we first introduce
the definition of weak solution to the problem (1.2).
\begin{definition1} A function $u(t,x) \in C(\mathbb{R}^{+}\times \mathbb{S})\cap L^{\infty}(\mathbb{R}^{+}; H^{1}(\mathbb{S}))$ is said to be an
admissible global weak solution to (1.2)  if
 $u$ satisfies the equations in (1.2) and
$z(t,\cdot)\rightarrow z_0$ as $t\rightarrow 0^+$ in the sense of
distributions on $\mathbb{R}_+\times\mathbb{R}$. Moreover,
$$\mu(u)=\mu(u_{0})e^{-\lambda t}\ \ \ \ \ \text{and} \ \ \ \ \  \|u_x(t,\cdot)\|_{L^2(\mathbb{S})}=e^{-\lambda t}
\|u_{0,x}\|_{L^2(\mathbb{S})}.$$
\end{definition1}

The main result of this paper can be stated as follows:
\begin{theorem1} Let $u_{0}\in H^{1}(\mathbb{S}).$ Assume that $y_{0}=(\mu(u_{0})-u_{0,xx})\in M^{+}(\mathbb{S}),$ then the equation (1.2) has an admissible global weak solution in the sense of Definition 1.1. Moreover,
$$u\in L^{\infty}_{loc}(\mathbb{R}_+; W^{1,\infty}({\mathbb{S}}))\cap H_{loc}^{1}(\mathbb{R}^{+}\times \mathbb{S}).$$ Furthermore,
$y=(\mu(u)-u_{xx}(t,\cdot))\in M^{+}(\mathbb{S})$ for a.e. $t\in \mathbb{R}^{+}$ is uniformly bounded on $\mathbb{S}.$
 \end{theorem1}

\begin{remark1}
If $y_{0}=(\mu(u_{0})-u_{0,xx})\in M^{-}(\mathbb{S}),$ then the conclusions in Theorem 1.1 also hold with $y=(\mu(u)-u_{xx}(t,\cdot))\in M^{-}(\mathbb{S})$.
\end{remark1}

 The paper is organized as follows. In
 Section 2, we recall some useful lemmas and derive some priori estimates on global strong solutions to (1.2). In
Section 3, we obtain the global existence of approximate solutions
to (1.2) with smooth approximate initial data. In Section 4,
using Helly$^{,}$s theorem, we prove the
existence of global weak solutions to (1.2). Moreover, complete the proof of Theorem 1.1.

\section{Preliminaries}
\newtheorem {remark2}{Remark}[section]
\newtheorem{theorem2}{Theorem}[section]
\newtheorem{lemma2}{Lemma}[section]

On one hand, with $y=\mu(u)-u_{xx},$ the first equation in (1.2) takes the form of a quasi-linear evolution
equation of hyperbolic type:
\begin{equation}
u_{t}+uu_{x}=-\partial_{x}A^{-1}(2\mu(u)u+\frac{1}{2}u_{x}^{2})-\lambda u,
\end{equation}
where $A=\mu-\partial_{x}^{2}$ is an isomorphism between
$H^{s}$ and $H^{s-2}$ with the inverse
$v=A^{-1}w$ given explicitly by \cite{escher, k-l-m}
\begin{align*}
v(x)&=(\frac{x^{2}}{2}-\frac{x}{2}+\frac{13}{12})\mu(w)+(x-\frac{1}{2})\int_{0}^{1}\int_{0}^{y}w(s)dsdy\\
\nonumber&-\int_{0}^{x}\int_{0}^{y}w(s)dsdy+\int_{0}^{1}\int_{0}^{y}\int_{0}^{s}w(r)drdsdy.
\end{align*}
Since $A^{-1}$ and
$\partial_{x}$ commute, the following identities hold
$$
A^{-1}\partial_{x}w(x)=(x-\frac{1}{2})\int_{0}^{1}w(x)dx-\int_{0}^{x}w(y)dy+\int_{0}^{1}\int_{0}^{x}w(y)dydx,
$$
and
\begin{equation}
A^{-1}\partial_{x}^{2}w(x)=-w(x)+\int_{0}^{1}w(x)dx.
\end{equation}
On the other hand, integrating both sides of the first equation in (1.2) with respect
to $x$ on $\mathbb{S}$, we obtain
$$\frac{d}{dt}\mu(u)=-\lambda \mu(u),$$ it follows that
\begin{equation}
\mu(u)=\mu(u_{0})e^{-\lambda t}:=\mu_{0}e^{-\lambda t},
\end{equation}
where
$$\mu_{0}:=\mu(u_{0})=\int_{\mathbb{S}}u_{0}(x)dx.$$
Combining  (2.1) and (2.3), the equation (1.2) can be rewrite as
\begin{equation}
\left\{\begin{array}{ll}
u_{t}+uu_{x}=-\partial_{x}A^{-1}(2\mu_{0}e^{-\lambda t}u+\frac{1}{2}u_{x}^{2})-\lambda u,&t > 0,\,x\in \mathbb{R},\\
u(0,x) = u_{0}(x),& x\in \mathbb{R}, \\
u(t,x+1)=u(t,x), & t \geq 0, x\in \mathbb{R}.\\\end{array}\right. \\
\end{equation}
If we rewrite the inverse of the operator $A=\mu-\partial_{x}^{2}$
in terms of a Green's function, we find
$(A^{-1}m)(x)=\int_{0}^{1}g(x-x^{'})m(x^{'})dx^{'}=(g*m)(x).$ So, we
get another equivalent form:
\begin{equation}
\left\{\begin{array}{ll}
u_{t}+uu_{x}=-\partial_{x}g*(2\mu_{0}e^{-\lambda t}u+\frac{1}{2}u_{x}^{2})-\lambda u,&t > 0,\,x\in \mathbb{R},\\
u(0,x) = u_{0}(x),& x\in \mathbb{R}, \\
u(t,x+1)=u(t,x), & t \geq 0, x\in \mathbb{R}.\\\end{array}\right. \\
\end{equation}
where the Green's function $g(x)$ is given \cite{l-m-t} by
\begin{equation}
g(x)=\frac{1}{2}x(x-1)+\frac{13}{12} \ \ \ \text{for} \
x\in[0,1)\simeq \mathbb{S},
\end{equation}
and is extended periodically to the real line. In other words,
$$g(x-x^{'})=\frac{(x-x^{'})^{2}}{2}-\frac{|x-x^{'}|}{2}+\frac{13}{12},
\ \ \ x,x^{'}\in[0,1)\simeq \mathbb{S}.$$ In particular, $\mu(g)=1.$

Given $u_{0}\in H^{s}$ with $s>\frac{3}{2}$, Theorem 2.2 in \cite{l-y}
ensures the existence of a maximal  $T> 0$ and a solution $u$ to (1.2) such that
$$u=u(\cdot,u_{0})\in C([0,T); H^{s}(\mathbb{S}))\cap C^{1}([0,T);H^{s-1}(\mathbb{S})).$$
Consider now the following initial value problem
\begin{equation}
\left\{\begin{array}{ll}q_{t}=u(t,q),\ \ \ \ t\in[0,T), \\
q(0,x)=x,\ \ \ \ x\in\mathbb{R}. \end{array}\right.
\end{equation}

\begin{lemma2}\cite{l-y}
Let $u_{0}\in H^{s}$ with $s>\frac{3}{2},$ $T> 0$ be the maximal existence time. Then Eq.(2.7) has a unique solution
$q\in C^1([0,T)\times \mathbb{R};\mathbb{R})$ and the map
$q(t,\cdot)$ is an increasing diffeomorphism of $\mathbb{R}$ with
$$
q_{x}(t,x)=exp\left(\int_{0}^{t}u_{x}(s,q(s,x))ds\right)>0, \ \
(t,x)\in [0,T)\times \mathbb{R}.$$
Moreover, with $y=\mu(u)-u_{xx},$ we have
$$
y(t,q(t,x))q_{x}^{2}=y_{0}(x)e^{-\lambda t}.
$$
\end{lemma2}

\begin{lemma2} If $y_{0}=\mu_{0}-u_{0,xx}\in H^{1}(\mathbb{S})$ does not change sign, then the corresponding solution $u$ to (2.5) of the initial value $u_{0}$ exists globally in time, that is $u\in C(\mathbb{R}^{+}, H^{3}(\mathbb{S}))\cap C^{1}(\mathbb{R}^{+}, H^{2}(\mathbb{S})).$ Moreover, the following properties hold:\\
(1) $\mu(u)=\mu_{0}e^{-\lambda t}, \ \ t\in[0, \infty),$\\
(2) $\|u_{x}\|_{L^{2}(\mathbb{S})}^{2}=e^{-2\lambda t}\mu_{1}^{2},\ \ t \in[0, \infty), \ \ \text{with} \ \ \mu_{1}=\left(\int_{\mathbb{S}}u_{0,x}^{2}dx\right)^{\frac{1}{2}},$\\
(3) $\|u(t,\cdot)\|_{L^{\infty}(\mathbb{S})}\leq|\mu_{0}|+\frac{\sqrt{3}}{6}\mu_{1},$\\
(4) $y(t,x), \ u(t,x)$ have the same sign with $y_{0}(x),$ and $\|u_{x}\|_{L^{\infty}(\mathbb{R}^{+}\times \mathbb{S})}\leq |\mu_{0}|,$\\
(5) $|\mu_{0}|e^{-\lambda t}=\|y_{0}\|_{L^{1}(\mathbb{S})}e^{-\lambda t}=\|y(t,\cdot)\|_{L^{1}(\mathbb{S})}=\|u(t,\cdot)\|_{L^{1}(\mathbb{S})}.$
\end{lemma2}
\textbf{Proof} Except for (4) and (5), all of the conclusions in Lemma 2 can be found in \cite{l-y}. So we only need to prove (4) and (5) here. Firstly, Lemma 2.1 and $u=g*y,$ $g\geq 0$ imply $y(t,x), \ u(t,x)$ have the same sign with $y_{0}(x).$ Moreover, from the proof of Theorem 5.1 in \cite{l-y}, we have $u_{x}(t,x)\geq - |\mu_{0}|.$ Now note that given $t\in[0,T),$ there is a $\xi(t)\in \mathbb{S}$ such that $u_{x}(t, \xi(t))=0$ by the periodicity of $u$ to $x$-variable. If $y_{0}\geq 0,$ then $y\geq 0.$ For $x\in[\xi(t), \xi(t)+1],$ we have
\begin{align*}
u_{x}(t,x)=\int_{\xi(t)}^{x}\partial_{x}^{2}u(t,x)dx=&\int_{\xi(t)}^{x}(\mu(u)-y)dx=\mu(u)(x-\xi(t))-\int_{\xi(t)}^{x}ydx\\
\leq \ &\mu(u)(x-\xi(t))\leq |\mu_{0}|.
\end{align*}
It follows that $u_{x}(t,x)\leq|\mu_{0}|.$
 On the other hand, if $y_{0}\leq 0,$ then $y\leq 0.$ Therefore, for $x\in[\xi(t), \xi(t)+1],$ we have
\begin{align*}
u_{x}(t,x)=\int_{\xi(t)}^{x}\partial_{x}^{2}u(t,x)dx=&\int_{\xi(t)}^{x}(\mu(u)-y)dx\leq\mu(u)(x-\xi(t))-\int_{\mathbb{S}}ydx\\
= \ &\mu(u)(x-\xi(t))-\int_{\mathbb{S}}(\mu(u)-u_{xx})dx=\mu(u)(x-\xi(t)-1)\leq |\mu_{0}|.
\end{align*}
It follows that $u_{x}(t,x)\leq |\mu_{0}|.$ So we have $\|u_{x}\|_{L^{\infty}(\mathbb{R}^{+}\times \mathbb{S})}\leq |\mu_{0}|,$ this complete the proof of (4). By the first equation of (1.1), we have
$$\int_{\mathbb{S}}y(t,x)dx=(\int_{\mathbb{S}}y_{0}(x)dx)e^{-\lambda t}=\mu_{0}e^{-\lambda t}.$$
If $y_{0}\geq 0,$ then $y\geq 0$ and $\mu_{0}\geq 0,$ we have
$$\|y\|_{L^{1}(\mathbb{S})}=\int_{\mathbb{S}}y(t,x)dx=(\int_{\mathbb{S}}y_{0}(x)dx)e^{-\lambda t}=\|y_{0}\|_{L^{1}(\mathbb{S})}e^{-\lambda t}=\mu_{0}e^{-\lambda t}.$$
If $y_{0}\leq 0,$ then $y\leq 0$ and $\mu_{0}\leq 0,$ we have
$$\|y\|_{L^{1}(\mathbb{S})}=-\int_{\mathbb{S}}y(t,x)dx=(\int_{\mathbb{S}}(-y_{0}(x))dx)e^{-\lambda t}=\|y_{0}\|_{L^{1}(\mathbb{S})}e^{-\lambda t}=-\mu_{0}e^{-\lambda t}.$$ Combining this two equalities above, we have $|\mu_{0}|e^{-\lambda t}=\|y_{0}\|_{L^{1}(\mathbb{S})}e^{-\lambda t}=\|y(t,\cdot)\|_{L^{1}(\mathbb{S})}.$ A similar discuss implies $\|y(t,\cdot)\|_{L^{1}(\mathbb{S})}=\|u(t,\cdot)\|_{L^{1}(\mathbb{S})}.$ This complete the proof of (5). $\Box$

\begin{lemma2}(\cite{s}) Assume $X\subset B \subset Y$ with compact
imbedding $X\rightarrow B$ ($X,B$ and $Y$ are Banach spaces), $1\leq
p\leq \infty$ and (1) $F$ is bounded in $L^{p}(0,T; X)$, (2)
$\|\tau_{h}f-f\|_{L^{p}(0,T-h; Y)}\rightarrow 0$ as $h\rightarrow 0$
uniformly for $f\in F.$ Then $F$ is relatively compact in
$L^{p}(0,T; B)$ (and in $C(0,T; B)$ if $p=\infty$), where
$(\tau_{h}f)(t)=f(t+h)$ for $h>0,$ if $f$ is defined on $[0,T],$
then the translated function $\tau_{h}f$ is defined on $[-h, T-h].$
\end{lemma2}

\begin{lemma2}(Helly$^{\prime}$s theorem \cite{ip})
Let an infinite family of function $F={f(x)}$ be defined on the segment $[a,b].$ If all functions of the family and the total variation of all functions of the family are bounded by a single number $|f(x)|\leq K, \ \ \ \bigvee\limits_{a}^{b}(f)\leq K,$ then there exists a sequence ${f_{n}(x)}$ in the family $F$ which converges at very point of $[a,b]$ to some function $\varphi(x)$ of finite variation.
\end{lemma2}

\begin{lemma2}(\cite{mnrr}) Let $T>0.$ If $f,g\in L^{2}((0,T);H^{1}(\mathbb{R}))$ and $\frac{df}{dt}, \frac{dg}{dt}\in  L^{2}((0,T);H^{-1}(\mathbb{R})),$ then $f,g$ are a.e. equal to a function continuous from $[0,T]$ into $L^{2}(\mathbb{R})$ and
$$\langle f(t), g(t)\rangle-\langle f(s), g(s)\rangle=\int_{s}^{t}\left\langle \frac{d(f(\tau))}{d\tau}, g(\tau)\right\rangle d\tau+\int_{s}^{t}\left\langle \frac{d(g(\tau))}{d\tau}, f(\tau)\right\rangle d\tau$$ for all $s,t\in [0,T].$
\end{lemma2}

\section{Global approximate solutions}
\newtheorem {remark3}{Remark}[section]
\newtheorem{theorem3}{Theorem}[section]
\newtheorem{lemma3}{Lemma}[section]
\newtheorem{definition3}{Definition}[section]
\newtheorem{claim3}{Claim}[section]

In the section, we will prove the existence of global approximate
solutions and give some useful estimates to the approximate solutions. Now we consider the approximate equation of (2.5) as follows:

\begin{equation}
\left\{\begin{array}{ll}
u^{n}_{t}+u^{n}u^{n}_{x}=-\partial_{x}g*(2\mu^{n}_{0}e^{-\lambda t}u^{n}+\frac{1}{2}(u^{n}_{x})^{2})-\lambda u^{n},&t > 0,\,x\in \mathbb{R},\\
u^{n}(0,x) = u^{n}_{0}(x),& x\in \mathbb{R}, \\
u^{n}(t,x+1)=u^{n}(t,x), & t \geq 0, x\in \mathbb{R},\\\end{array}\right. \\
\end{equation}
where $u_{0}^{n}(x)=\phi_n\ast u_{0}\in H^{\infty}(\mathbb{S})$ for $n\geq 1$ and $\mu^{n}_{0}=\int_{\mathbb{S}}u_{0}^{n}(x)dx.$ Here $\{\phi_n\}_{n\geq 1}$ are the mollifiers
$$ \phi_n(x):=\left(\int_{\mathbb{R}}\phi(\xi)d\xi\right)^{-1}n\phi(nx),\ \ \ \ \ x\in\mathbb{R}, \ \ n\geq 1,$$
where $\phi\in C_c^{\infty}(\mathbb{R})$ is defined by
$$\phi(x)=\left \{\begin {array}{ll}e^{1/(x^2-1)},  \ \ \ \ \ \ \ \ |x|<1,
 \\ 0, \ \ \ \ \ \ \ \ \ \ \ \ \ \ \ \ \ \ |x|\geq 1.\end {array}\right.$$ Obviously, $\|\phi_n\|_{L^{1}(\mathbb{R})}=1.$ Clearly, we have
\begin{align}
 u_0^n\rightarrow u_0 \ \ \text{in}  \ \
H^1(\mathbb{S}), \ \ \ \text{as}\ \ n\rightarrow\infty
\end{align} and
\begin{align}
&\|u_0^n\|_{L^2(\mathbb{S})}\leq \|u_0\|_{ L^2(\mathbb{S})}, \ \ \ \ \ \ \
\|u_{0,x}^n\|_{L^2(\mathbb{S})}\leq \|u_{0,x}\|_{L^2(\mathbb{S})},\\
\nonumber & \|u_0^n\|_{H^{1}(\mathbb{S})}\leq \|u_0\|_{H^{1}(\mathbb{S})},\ \ \ \ \ \ \ \|u_0^n\|_{L^1(\mathbb{S})}\leq \|u_0\|_{ L^1(\mathbb{S})}
\end{align}
in view of Young$^{,}$s inequality. Note that
\begin{align*}
\mu_{0}^{n}=\mu(u_{0}^{n})=\int_{\mathbb{S}}u_{0}^{n}(x)dx=&\int_{\mathbb{S}}\int_{\mathbb{R}}\phi_n(y)u_{0}(x-y)dydx
=\int_{\mathbb{R}}\int_{\mathbb{S}}\phi_n(y)u_{0}(x-y)dxdy\\
= &\int_{\mathbb{R}}\phi_n(y)\cdot\left(\int_{\mathbb{S}}u_{0}(x-y)dx\right)dy\\
=&\int_{\mathbb{R}}\phi_n(y)\cdot\left(\int_{\mathbb{S}}u_{0}(z)dz\right)dy\\
=&\int_{\mathbb{R}}\phi_n(y)\mu(u_{0})(x-y)dy\\
=&\phi_n*\mu(u_{0})=\mu(u_{0})=\mu_{0}.
\end{align*}
We can rewrite (3.1) as follows:
\begin{equation}
\left\{\begin{array}{ll}
u^{n}_{t}+u^{n}u^{n}_{x}=-\partial_{x}g*(2\mu_{0}e^{-\lambda t}u^{n}+\frac{1}{2}(u^{n}_{x})^{2})-\lambda u^{n},&t > 0,\,x\in \mathbb{R},\\
u^{n}(0,x) = u^{n}_{0}(x),& x\in \mathbb{R}, \\
u^{n}(t,x+1)=u^{n}(t,x), & t \geq 0, x\in \mathbb{R}.\\\end{array}\right. \\
\end{equation}
Moreover, for all $n\geq1,$ $y_{0}^{n}=\mu(u_{0}^{n})-u^{n}_{0,xx}=\mu_{0}-u^{n}_{0,xx}\in H^{1}(\mathbb{S})$ and
$$y_{0}^{n}=\mu(u^{n}_{0})-u^{n}_{0,xx}=\phi_n*\mu(u_{0})-\phi_n*u_{0,xx}=\phi_n*y_{0}\geq 0.$$
Thus, by Lemma 2.2, we obtain the corresponding solution $u^n\in
 C(\mathbb{R}^{+};H^3(\mathbb{S}))\cap C^1(\mathbb{R}^{+};H^{2}(\mathbb{S}))$ to Eq.(3.4) with initial data $u^n_0(x)$
and $y^{n}=\mu(u^{n})-u^{n}_{xx}\geq 0$, $u^{n}=g*y^{n}\geq 0$ for all $(t,x)\in \mathbb{R}^{+}\times \mathbb{S}).$ Furthermore, combining Lemma 2.2 and (3.3), we have:
\begin{align}
&\mu(u^{n})=\mu^{n}_{0}e^{-\lambda t}=\mu_{0}e^{-\lambda t}, \ \ t\in[0, \infty),\\
&\|u^{n}_{x}\|_{L^{2}(\mathbb{S})}^{2}=e^{-2\lambda t}\|u^{n}_{0,x}\|_{L^{2}(\mathbb{S})}^{2}\leq\|u_{0,x}\|_{L^{2}(\mathbb{S})}^{2}=\mu_{1}^{2},\ \ t \in[0, \infty),\\
&\|u^{n}(t,\cdot)\|_{L^{\infty}(\mathbb{S})}\leq|\mu^{n}_{0}|+\frac{\sqrt{3}}{6}\|u^{n}_{0,x}\|_{L^{2}(\mathbb{S})}\leq |\mu_{0}|+\frac{\sqrt{3}}{6}\|u_{0,x}\|_{L^{2}(\mathbb{S})}=|\mu_{0}|+\frac{\sqrt{3}}{6}\mu_{1},\\
&\|u^{n}_{x}\|_{L^{\infty}(\mathbb{R}^{+}\times \mathbb{S})}\leq |\mu^{n}_{0}|=|\mu_{0}|,\\
&|\mu_{0}|e^{-\lambda t}=\|y^{n}(t,\cdot)\|_{L^{1}(\mathbb{S})}=\|u^{n}(t,\cdot)\|_{L^{1}(\mathbb{S})}.
\end{align}

\section{Proof of Theorem 1.1}
\newtheorem{theorem4}{Theorem}[section]
\newtheorem{lemma4}{Lemma}[section]
\newtheorem {remark4}{Remark}[section]
\newtheorem{corollary4}{Corollary}[section]

In this section, with the basic energy estimate in Section 3, we are
ready to obtain the necessary compactness of the approximate
solutions $u^n(t,x).$ Acquiring the precompactness of approximate
solutions, we prove the existence of the global weak solutions to the
equation (1.1).
\begin{lemma4}
For any fixed $T>0$, there exist a subsequence $\{u^{n_k}(t,x)\}$ of
the sequence $\{u^{n}(t,x)\}$ and some function  $u(t,x)\in
L^\infty(\mathbb{R}^{+}; H^1(\mathbb{S}))\cap H^{1}([0,T]\times \mathbb{S})$ such that
\begin{equation} u^{n_k}\rightharpoonup u \ \ \ \ \text{in}\
H^1([0,T]\times\mathbb{S}) \ \text{as}\ n_k\rightarrow\infty, \forall T>0,
\end{equation}
and
\begin{equation} u^{n_k}\rightarrow u \ \ \text{in} \ \
L^\infty([0,T]\times\mathbb{S}) \ \
 \text{as} \ \ n_k\rightarrow\infty.
\end{equation}
Moreover, $u(t,x)\in C(\mathbb{R}^{+}\times\mathbb{S}).$
\end{lemma4}
\textbf{Proof} Firstly, we will prove that the sequence $\{u^{n}(t,x)\}$ is uniformly
bounded in the space $H^{1}([0,T]\times \mathbb{S}).$ By (3.6)-(3.7), we have
\begin{equation}
\|u^{n}\|^{2}_{L^{2}([0,T]\times\mathbb{S})}=\int_{0}^{T}\int_{\mathbb{S}}(u^{n})^{2}dxdt=\int_{0}^{T} \|u^{n}\|^{2}_{L^{2}(\mathbb{S})}dx\leq \left(|\mu_{0}|+\frac{\sqrt{3}}{6}\mu_{1}\right)^{2}T
\end{equation}
\begin{equation}
\|u_{x}^{n}\|^{2}_{L^{2}([0,T]\times\mathbb{S})}=\int_{0}^{T}\int_{\mathbb{S}}(u_{x}^{n})^{2}dxdt=\int_{0}^{T} \|u_{x}^{n}\|^{2}_{L^{2}(\mathbb{S})}dx\leq \mu_{1}^{2}T
\end{equation}
Moreover, by (3.7)-(3.8), we obtain
\begin{equation}
\|u^{n}u^{n}_{x}\|_{L^{2}([0,T]\times \mathbb{S})}\leq
 \|u^{n}\|_{L^{2}([0,T]\times \mathbb{S})}\|u^{n}_{x}\|_{L^{\infty}([0,T]\times \mathbb{S})}\leq \left(|\mu_{0}|+\frac{\sqrt{3}}{6}\mu_{1}\right)|\mu_{0}|
\end{equation}
\begin{align}
\nonumber &\|\partial_{x}g*(2\mu_{0}e^{-\lambda t}u^{n}+\frac{1}{2}(u_{x}^{n})^{2})\|_{L^{2}([0,T]\times
\mathbb{S})}\\
\nonumber\leq \ & \|\partial_{x}g\|_{L^{2}([0,T]\times
\mathbb{S})}\|2\mu_{0}e^{-\lambda t}u^{n}+\frac{1}{2}(u_{x}^{n})^{2}\|_{L^{1}([0,T]\times
\mathbb{S})}\\
\nonumber\leq \ &
\frac{T}{12}\int_{0}^{T}\int_{\mathbb{S}}\left(2|\mu_{0}||u^{n}|+\frac{1}{2}(u_{x}^{n})^{2}\right)dxdt\\
 \leq \ &
\frac{T^{2}}{12}\left[\mu_{0}^{2}+(|\mu_{0}|+\frac{\sqrt{3}}{6}\mu_{1})^{2}+\mu_{1}^{2}\right].
\end{align}
Combining (4.3), (4.5)-(4.6) with Eq.(3.4), we know that
$\{u^{n}_{t}(t,x)\}$ is uniformly bounded in $L^{2}([0,T]\times
\mathbb{S}).$ Thus, (4.3), (4.4) and this conclusion imply that
$$\int_{0}^{T}\int_{\mathbb{S}}((u^{n})^{2}+(u^{n}_{x})^{2}+(u^{n}_{t})^{2})dxdt \leq K,$$
where $K=K(|\mu_{0}|, \mu_{1}, T, \lambda)\geq 0.$ It follows that $\{u^{n}(t,x)\}$ is uniformly
bounded in the space $H^{1}([0,T]\times \mathbb{S}).$ Thus (4.1) holds for some $u\in H^{1}([0,T]\times
\mathbb{S}).$

Observe that, for each $0\leq s, t\leq T,$
$$\|u^{n}(t,\cdot)-u^{n}(s,\cdot)\|_{L^{2}(\mathbb{S})}^{2}
=\int_{\mathbb{S}}(\int_{s}^{t}\frac{\partial u^{n}}{\partial
\tau}(\tau, x)d\tau)^{2}dx\leq
|t-s|\int_{0}^{T}\int_{\mathbb{S}}(u_{t}^{n})^{2}dxdt.$$ Note that $\{u^{n}(t,x)\}$ is uniformly bounded in $L^{\infty}([0,T];
H^{1}(\mathbb{S})),$ $\{u^{n}_{t}(t,x)\}$ is uniformly bounded in $L^{2}([0,T]\times
\mathbb{S})$ and $H^{1}(\mathbb{S})\subset\subset C(\mathbb{S})\subset
L^{\infty}(\mathbb{S})\subset L^{2}(\mathbb{S}),$ then (4.2) and
$u(t,x)\in C(\mathbb{R}^{+}\times\mathbb{S})$ is consequence of Lemma
2.3. $\Box$

Consequently, we will deal with $\partial_{x}g*(2\mu_{0}e^{-\lambda t}u^{n}+\frac{1}{2}(u^{n}_{x})^{2}).$ By (3.5), (3.7)-(3.9), we have that for fixed $t\in [0,T]$ the sequence $u^{n_{k}}_{x}(t,\cdot)\in BV(\mathbb{S})$ satisfies
$$\mathbb{V}(u^{n_{k}}_{x}(t,\cdot))=\|u^{n_{k}}_{xx}(t,\cdot)\|_{L^{1}(\mathbb{S})}=\|\mu(u^{n_{k}})-y^{n_{k}}\|_{L^{1}(\mathbb{S})}\leq \|\mu(u^{n_{k}})\|_{L^{1}(\mathbb{S})}+\|u^{n_{k}}\|_{L^{1}(\mathbb{S})}\leq 2|\mu_{0}|+\frac{\sqrt{3}}{6}\mu_{1}$$
and
$$\|u^{n_{k}}_{x}(t,\cdot)\|_{L^{\infty}(\mathbb{S})}\leq |\mu_{0}|\leq 2|\mu_{0}|+\frac{\sqrt{3}}{6}\mu_{1}.$$
Applying Lemma 2.4, we obtain that there exists a subsequence, denoted
again $\{u^{n_{k}}_{x}(t,\cdot)\},$ which converges at every point to some function $v(t,x)$ of finite
variation with $\mathbb{V}(v(t,\cdot))\leq 2|\mu_{0}|+\frac{\sqrt{3}}{6}\mu_{1}.$ Since for almost all $t\in [0,T],$ $u^{n_{k}}_{x}(t,\cdot)\rightarrow u_{x}(t,\cdot)$ in $\mathcal{D}^{\prime}(\mathbb{S})$ in view of Lemma 4.1, it follows that $v(t,\cdot)=u_{x}(t,\cdot)$ for a.e.$t\in [0, T].$ So we have \begin{equation}
u^{n_{k}}_{x}(t,\cdot)\rightarrow u_{x}(t,\cdot) \  a.e. \ \text{on} \ [0,T]\times\mathbb{S}, \ \ \ \text{for}  \ \  n_{k}\rightarrow \infty,
\end{equation}
and for a.e. $t\in [0,T],$
\begin{equation}
\mathbb{V}(u_{x}(t,\cdot))=\|u_{xx}(t,\cdot)\|_{M(\mathbb{S})}\leq 2|\mu_{0}|+\frac{\sqrt{3}}{6}\mu_{1}.
\end{equation}
Therefor,
\begin{align*}
&\|\partial_{x}g*(2\mu_{0}e^{-\lambda t}u^{n_{k}}+\frac{1}{2}(u^{n_{k}}_{x})^{2})-\partial_{x}g*(2\mu_{0}e^{-\lambda t}u+\frac{1}{2}(u_{x})^{2})\|_{L^{\infty}([0,T]\times \mathbb{S})}\\
\leq \ & \|\partial_{x}g\|_{L^{1}([0,T]\times \mathbb{S})}\|2\mu_{0}e^{-\lambda t}(u^{n_{k}}-u)+\frac{1}{2}(u^{n_{k}}_{x})^{2}-(u_{x})^{2}))\|_{L^{\infty}([0,T]\times \mathbb{S})}\\
\leq \ & \frac{T}{4}\left(2|\mu_{0}|\|u^{n_{k}}-u\|_{L^{\infty}([0,T]\times \mathbb{S})}+\frac{1}{2}\|u_{x}^{n_{k}}+u_{x}\|_{L^{\infty}([0,T]\times \mathbb{S})}\|u_{x}^{n_{k}}-u_{x}\|_{L^{\infty}([0,T]\times \mathbb{S})}\right).
\end{align*}
Combining this inequality with (4.2), (4.7) and note that
$$\|u_{x}(t,\cdot)\|_{L^{\infty}(\mathbb{S})} \leq \lim\limits_{n_{k}\rightarrow \infty} \|u^{n_{k}}_{x}(t,\cdot)\|_{L^{\infty}(\mathbb{S})}\leq |\mu_{0}|,$$
we have \begin{equation}
\partial_{x}g*(2\mu_{0}e^{-\lambda t}u^{n_{k}}+\frac{1}{2}(u^{n_{k}}_{x})^{2})\rightarrow \partial_{x}g*(2\mu_{0}e^{-\lambda t}u+\frac{1}{2}u_{x}^{2})
\end{equation}
a.e. on $[0,T]\times \mathbb{S}.$ The relations (4.2), (4.7) and (4.9) imply that $u$ satisfies Eq.(2.5) in $\mathcal{D}^{\prime}([0,T]\times \mathbb{S}).$ Moreover, since $$\|u(t,\cdot)\|_{L^{\infty}(\mathbb{S})} \leq \lim\limits_{n_{k}\rightarrow \infty} \|u^{n_{k}}(t,\cdot)\|_{L^{\infty}(\mathbb{S})}\leq |\mu_{0}|+\frac{\sqrt{3}}{6}\mu_{1},$$ we have $u\in L_{loc}^{\infty}(\mathbb{R}^{+}, W^{1,\infty}(\mathbb{S}))$ in view of $T$ in (4.2) and (4.7) being arbitrary.

Now, we prove that $\mu(u)=\mu(u_{0})e^{-\lambda t},$ $\|u_{x}\|_{L^{2}(\mathbb{S})}^{2}=e^{-2\lambda t}\|u_{0,x}\|_{L^{2}(\mathbb{S})}^{2}$ and $(\mu(u)-u_{xx}(t,\cdot))\in M^{+}(\mathbb{S})$ is uniformly bounded on $\mathbb{S}.$

On one hand, by (4.2), we have$$\int_{\mathbb{S}}u^{n_{k}}(t,x)dx\rightarrow \int_{\mathbb{S}}u(t,x)dx=\mu(u) \ \ \ \text{as} \ n_{k}\rightarrow \infty.$$
On the other hand, $$\int_{\mathbb{S}}u^{n_{k}}(t,x)dx=\mu(u^{n_{k}})=\mu_{0}e^{-\lambda t}.$$ Obviously, $\mu(u)=\mu(u_{0})e^{-\lambda t}$ by the uniqueness of limit.

By $u$ satisfies (2.5) in the sense of distribution, we have
\begin{equation}
\phi_n\ast u_{t}+\phi_n\ast(uu_{x})=-\phi_n\ast\left(\partial_{x}g*(2\mu(u)u+\frac{1}{2}u_{x}^{2})\right)-\lambda \phi_n\ast u.
\end{equation}
Differentiating (4.10) with respect to $x,$ we have
$$(\phi_n\ast u_{x})_{t}+\phi_n\ast(uu_{xx})=-\phi_n\ast\left(2(\mu(u))^{2}+\frac{1}{2}\mu(u_{x}^{2})-2\mu(u)u+\frac{1}{2}u_{x}^{2}\right)-\lambda \phi_n\ast u_{x},$$ here we used the formula (2.2).
Multiplying the equality above with $\phi_n\ast u_{x}$ and integrating the result with respect to $x$ on $\mathbb{S},$ we obtain
\begin{align*}
&\frac{1}{2}\frac{d}{dt}\int_{\mathbb{S}}(\phi_n\ast u_{x})^{2}dx+\int_{\mathbb{S}}(\phi_n\ast u_{x})(\phi_n\ast(uu_{xx}))dx\\
=&-\int_{\mathbb{S}}(\phi_n\ast u_{x})(\phi_n\ast(2(\mu(u))^{2}+\frac{1}{2}\mu(u_{x}^{2})-2\mu(u)u+\frac{1}{2}u_{x}^{2}))dx-\lambda\int_{\mathbb{S}}(\phi_n\ast u_{x})^{2}dx.
\end{align*}
Note that
$$\int_{\mathbb{S}}(\phi_n\ast u_{x})(\phi_n\ast(2(\mu(u))^{2}+\frac{1}{2}\mu(u_{x}^{2}))dx=0$$
and
$$\int_{\mathbb{S}}(\phi_n\ast u_{x})(\phi_n\ast(-2\mu(u)u))dx=-2\mu(u)\int_{\mathbb{S}}(\phi_n\ast u_{x})(\phi_n\ast u)dx=0$$
we have
$$\frac{d}{dt}\int_{\mathbb{S}}(\phi_n\ast u_{x})^{2}dx+2\lambda\int_{\mathbb{S}}(\phi_n\ast u_{x})^{2}dx=-2\int_{\mathbb{S}}(\phi_n\ast u_{x})(\phi_n\ast(uu_{xx}))dx-\int_{\mathbb{S}}(\phi_n\ast u_{x})(\phi_n\ast u_{x}^{2})dx.$$
Let $$f_{n}(t)=\int_{\mathbb{S}}(\phi_n\ast u_{x})^{2}dx,$$ $$g_{n}(t)=-2\int_{\mathbb{S}}(\phi_n\ast u_{x})(\phi_n\ast(uu_{xx}))dx-\int_{\mathbb{S}}(\phi_n\ast u_{x})(\phi_n\ast u_{x}^{2})dx,$$ then we obtain for a.e. $t\in \mathbb{R}^{+},$
\begin{equation}
\frac{d f_{n}(t)}{dt}+2\lambda f_{n}(t)=g_{n}(t).
\end{equation}
Applying Lemma 2.5 to $\phi_n\ast u_{x},$ it follows from (4.11) that
\begin{equation}
f_{n}(t)-e^{-2\lambda t}f_{n}(0)=\int_{0}^{t}e^{-2\lambda (t-s)}g_{n}(s)ds.
\end{equation}
Since $g_{n}(t)\rightarrow 0$ as $n\rightarrow \infty$ for a.e. $t\in\mathbb{R}^{+},$ we have for any $T>0,$ there exists a constant $K(T)>0,$ such that $|g_{n}(t)|\leq K(T), \ \ t\in[0,T], \ n\geq 1.$ An application of Lebesgue$^{\prime}$s dominated convergence theorem to (4.12), we obtain
$$\lim\limits_{n\rightarrow\infty}[ f_{n}(t)-e^{-2\lambda t}f_{n}(0)]=0.$$ That is for all $t\in \mathbb{R}^{+},$ we have
$\|u_{x}\|_{L^{2}(\mathbb{S})}^{2}=e^{-2\lambda t}\|u_{0,x}\|_{L^{2}(\mathbb{S})}^{2}.$

Note that $L^{1}(\mathbb{S})\subset M(\mathbb{S}).$ By (4.8) and $\mu(u)=\mu_{0}e^{-\lambda t},$ we have
$$\|\mu(u)-u_{xx}(t,\cdot)\|_{M(\mathbb{S})}\leq \|\mu(u)\|_{L^{1}(\mathbb{S})}+\|u_{xx}(t,\cdot)\|_{M(\mathbb{S})}\leq 3|\mu_{0}|+\frac{\sqrt{3}}{6}\mu_{1}.$$ It follows that for all $t\in \mathbb{R}^{+},$ $(\mu(u)-u_{xx}(t,\cdot))\in M(\mathbb{S})$ is uniformly bounded on $\mathbb{S}.$ For any fixed $T>0,$ in view of (4.2) and (4.7), we have for all $t\in [0,T],$
$$[\mu(u^{n_{k}})-u^{n_{k}}_{xx}(t,\cdot)]\rightarrow [\mu(u)-u_{xx}(t,\cdot)] \ \ \ \text{in} \ \  \mathcal{D}^{\prime}(\mathbb{S}) \ \ \text{for} \ \ n\rightarrow\infty.$$ Since $\mu(u^{n_{k}})-u^{n_{k}}_{xx}(t,\cdot)=y^{n_{k}}(t,\cdot)\geq 0$ for all $(t,x)\in \mathbb{R}^{+}\times\mathbb{S},$ we have $(\mu(u)-u_{xx}(t,\cdot))\in L_{loc}^{\infty}(\mathbb{R}^{+}, M^{+}(\mathbb{S})).$

\bigskip

\noindent\textbf{Acknowledgments} This work was partially supported by NNSFC (No. 10971235), RFDP (No. 200805580014), NCET (No. 08-0579) and the key project of Sun Yat-sen University.

\end{document}